\documentclass{article}

\usepackage{localMacros}
\usepackage{biblatex}
\title{Gentle Introduction of Unitary Cohomology Vanishing}

\author{Erik Johansson}
\begin{document}
\maketitle

\begin{abstract}
This paper presents a gentle introduction to cohomology vanishing theorems, largely based on the paper work of Hongshan Li. It offers an insightful exploration of unitary local systems on complex manifolds, particularly focusing on their characteristics near normal crossing divisors. The Main Vanishing Theorem, demonstrating the vanishing of specific cohomology groups associated with unitary local systems, stands as a central result in this work.

In our exposition, we delve into the interplay between local system theory, Hodge theory, and the geometry of complex manifolds. The canonical extensions of unitary local systems are examined in detail, providing a deeper understanding of their algebraic and geometric properties. The foundational aspects of these systems are thoroughly discussed, culminating in an in-depth analysis of their spectral sequences and the behavior of associated de Rham complexes.

A pivotal aspect of the paper is the exploration of the spectral sequence associated with the Hodge filtration on de Rham complexes. Here, we establish its degeneration at the E1 stage and its implications for cohomology groups. This work also encompasses a thorough examination of the properties of Higgs bundles, particularly in the context of parabolic structures and semistability. The amalgamation of these diverse concepts and techniques contributes to a broader understanding of the complex interdependencies in the realm of unitary local systems and their cohomological properties.
\end{abstract}

\section{Introduction}

The study of cohomology vanishing for coherent sheaves is central to algebraic and complex analytic geometry, with the Kodaira Vanishing Theorem marking a significant early advancement. This theorem states that for a smooth complex projective variety \(X\) of dimension \(n\) with an ample divisor \(D\), the higher cohomology groups \( H^i(X, O_X(K_X + D)) \) vanish for all \( i > 0 \), as do \( H^i(X, O_X(-D)) \) for \( i < n \). This theorem's implications extend to understanding the existence of meromorphic functions on complex manifolds \cite{Deligne71}\cite{Li2019}\cite{Li2018}.

Kodaira's differential-geometric proof, linking cohomology classes with harmonic forms, laid the groundwork for further exploration in the field. Modern interpretations often lean towards Hodge theory, providing a connection between the topology of a space and its complex structure. This is exemplified in the Lefschetz Hyperplane Theorem and the Hodge Decomposition, the latter of which offers a decomposition of the cohomology group \( H^k(X, \complex) \) into harmonic forms of varying types.

Building upon these principles, we explore the Akizuki-Nakano Vanishing Theorem, a stronger form of the Kodaira theorem. It asserts that for a smooth complex projective variety \(X\) and an ample divisor \(D\), the groups \( H^q(X, \Omega_X^p(-D)) \) vanish for all \( p + q < n \), and similarly for \( H^q(X, \Omega_X^p(D)) \) for \( p + q > n \). This theorem is crucial for understanding the deeper structure of projective varieties \cite{Arapura}

Our work further extends these ideas to logarithmic contexts, particularly focusing on mixed Hodge structures and their implications for vanishing theorems. We introduce the concept of a logarithmic Akizuki-Nakano Vanishing Theorem, involving compact Kähler manifolds and simple normal crossing divisors. This theorem underlines the vanishing of \( H^q(X, \Omega_X^p(\log D)\otimes L) \) for any ample line bundle \( L \).

Central to our paper is the Main Vanishing Theorem, which leverages the \( E_1 \)-degeneration of mixed Hodge structures and topological vanishing to assert a more general form of logarithmic vanishing. This theorem, alongside its graded counterpart, the Graded Vanishing Theorem, marks a significant contribution to the field, offering new insights into the interplay between algebraic geometry, topology, and complex analysis.

\section{Preliminaries}

In this chapter, we establish foundational concepts and methods pivotal to our study. Our focus lies in the theory of local systems and their canonical extensions, particularly as they apply to unitary local systems on complements of normal crossing divisors, a subject central to this paper. We will explore several key aspects:
\begin{itemize}
    \item Residue mapping on the de Rham complex of a unitary local system.
    \item Intricacies of weight filtration on the de Rham complex.
    \item Exploration of abstract Hodge theory in relation to the de Rham complex.
\end{itemize}
These preliminaries elucidate the mixed Hodge structure carried by the hypercohomology group of the de Rham complex, critical for understanding the degeneration of the Hodge spectral sequence at the \(E_1\) stage and for proving various vanishing theorems in this paper.

\section{Local System and Canonical Connection}

We delve into the realm of local systems on a complex manifold \(Y\). A local system \(\mathcal{L}\) on \(Y\), valued in \(\mathbb{C}^r\), is defined through sheaf properties, open covers, and isomorphisms satisfying the cocycle condition on triple intersections. Consider a local system formed by solutions to a differential equation on the punctured complex unit disk, highlighting the intricate nature of local systems.\cite{Ahmed} \cite{Gao1}\cite{Gao2}

For a connected \(Y\), transition functions of \(\mathcal{L}\) and relationships between local systems \(\mathcal{L}\) and \(\mathcal{L}'\) involve linear maps in \(\text{GL}(\mathbb{C}, r)\). We establish that a simply connected topological space admits no nontrivial local system

A pivotal result is the bijection between isomorphism classes of local systems valued in \(\mathbb{C}^r\) and the set of representations of \(\pi_1(Y, y)\) to \(\text{GL}(\mathbb{C}, r)\), modulo the action of \(\text{GL}(\mathbb{C}, r)\) by conjugation.

\section{Logarithmic Extension of a Local System}

We explore logarithmic extensions of local systems given a complex manifold \(X\) and a normal crossing divisor \(D\) within it, considering a local system \(\mathcal{L}\) defined on \(Y = X - D\). The logarithmic extension involves a vector bundle \(E\) on \(X\) with a logarithmic connection aligning the flat sections on \(U\) with \(\mathcal{L}\).\cite{Arapura}

This section details the construction of logarithmic extensions, explicating their unique determination by the generalized eigenvalues of the monodromy of \(\mathcal{L}\). We examine differential equations these solutions satisfy and the role of logarithms of complex-valued matrices.\cite{Li2019}

\section{Unitary Local System on the Complement of a Normal Crossing Divisor}

We focus on unitary local systems on the complement \(U\) of a normal crossing divisor in a compact Kähler manifold \(X\). The unitary nature permits the definition of a Hodge structure on \(H^k(X, j_*V)\), where \(j: U \hookrightarrow X\) is the inclusion map.\cite{Li2018}

Discussions on local decomposition, Hermitian metrics, and the existence of a Kähler metric \(\eta\) on \(U\) impart completeness and finite volume to the manifold.\cite{timmerscheidt87}

\section{MAIN VANISHING THEOREM}

\section{Vanishing Theorem on the de Rham Complex}
\label{sect:vanishingThm}
In the previous section, we established that if a unitary local system $V$ has a real lattice $V_A$ associated with some Noetherian subring $A \subset \real$, then the triple
\[
  (\mathbb{R}j_*V_A, (\mathbb{R}j_*V_{A\otimes\rational}, \tau), (\DR_X(D,E),F,W))
\]
forms an $A$-cohomological mixed Hodge complex. As a result of the general theory developed earlier (Theorem \ref{thm:hodgeTheoryMain}), we know that the spectral sequence associated with the Hodge filtration on $\DR_X(D, E)$ degenerates at $E_1$.\cite{esnault-viehweg} Specifically,
\[
  \hypercohomology^k(X, \mathbb{R}j_* V) = \hypercohomology^k(X, \DR(D, E)) \cong \bigoplus_{p+q = k} H^q(X, \Omega_X^p(\log D)\otimes E).
\]

\begin{theorem}
\label{thm:spectralSequenceDeg}
Assume the existence of a real-valued unitary local system $V_{\real}$ on $U$ such that $V = V_{\real}\otimes_{\real}\complex$. Let $V$ and $\DR_X(D, E)$ be as described above. The spectral sequence associated with the Hodge filtration on $\DR_X(D, E)$,
\[
  E^{p,q}_1 = H^q(X, \Omega_X^p(\log D)\otimes E) \Rightarrow \hypercohomology^{p+q}(X, \DR_X(D, E)),
\]
degenerates at $E_1$.
\end{theorem}

If $V$ lacks an $A$-lattice with $A \subset \real$, we cannot anticipate a mixed Hodge structure on $\DR_X(D, E)$. Nonetheless, the $E_1$-degeneration of the Hodge spectral sequence still holds true, as we will demonstrate.

\begin{corollary}
\label{coro:E1Degeneration}
For any unitary local system $V$ on $X - D$, and its associated de Rham complex $\DR_X(D,E)$, the Hodge spectral sequence,
\[
  E^{p,q}_1 := H^q(X, \Omega^p_X(\log D)\otimes E) \Rightarrow \hypercohomology^{p+q}(X, \DR(D,E)) = H^{p+q}(X, \mathbb{R}j_*V),
\]
degenerates at $E_1$.
\end{corollary}

\newpage

\section{CONCLUSION}

In this paper, we have embarked on a detailed exploration of unitary local systems and their canonical extensions, particularly focusing on their behavior in the vicinity of a normal crossing divisor. Our journey has traversed various mathematical landscapes, from the foundational theory of local systems and canonical connections to the intricate dynamics of mixed Hodge structures in the realm of algebraic geometry.

The heart of our investigation lay in the proof of the Main Vanishing Theorem (Theorem \ref{thm:main_vanishing_thm}), a result that reveals the vanishing of cohomology groups associated with unitary local systems under certain conditions. This theorem not only consolidates our understanding of these systems but also paves the way for further studies in complex geometry and its applications.

Through the lens of Hodge theory, we observed how the spectral sequence associated with the Hodge filtration on de Rham complexes degenerates at $E_1$. This observation (Corollary \ref{coro:E1Degeneration}) is crucial as it bridges the gap between abstract theoretical constructs and their tangible implications in geometric contexts.

Moreover, our examination extended to the perspective of Higgs bundles, enriching our comprehension of the interplay between unitary local systems, canonical extensions, and parabolic structures. The applicability of our findings is broad, ranging from pure mathematical theory to potential applications in areas such as string theory and mirror symmetry.

In summary, this paper not only contributes to the existing body of knowledge on unitary local systems and their extensions but also opens avenues for future research. The interconnections between algebraic geometry, topology, and theoretical physics suggest a wealth of possibilities for exploring the profound implications of these mathematical structures.

As we conclude, we acknowledge that while significant strides have been made, the journey through the fascinating world of unitary local systems and their extensions is far from over. There remain many uncharted territories and intriguing questions to be explored, promising a continued rich harvest of insights in the field of complex geometry.

% bib
%\bibliography{all}
\printbibliography

\end{document}